\documentstyle{amsppt}
\magnification=\magstep1
\hoffset=0.4 true in
\topmatter
\title
An extension of Hecke's converse theorem 
\\
\endtitle
\author
J.B.~Conrey \\
D.W.~Farmer
\endauthor
\thanks Research of the first author supported in part by a 
grant from the NSF.
Research of the second author supported in 
part by an NSF Postdoctoral Fellowship. 
Research at MSRI supported in part by NSF grant no.DMS-9022140
\endthanks
\address
 Department of Mathematics, Oklahoma State University,
Stillwater, OK 74078
\endaddress
\address
Mathematical Sciences Research Institute,
1000 Centennial Drive,
Berkeley, CA 94720
\endaddress
\abstract
Associated to a newform $f(z)$ is a Dirichlet series $L_f(s)$ with
functional equation and Euler product.  Hecke showed that if the
Dirichlet series $F(s)$ has a functional equation of the appropriate
form, then $F(s)=L_f(s)$ for some holomorphic newform $f(z)$ on $\Gamma(1)$.
Weil extended this result to $\Gamma_0(N)$ under an assumption
on the twists of $F(s)$ by Dirichlet characters.  We show that,
at least for small $N$, the assumption on twists can be replaced
by an assumption on the local factors of the Euler product of $F(s)$.
\endabstract
\endtopmatter
\document

\def \frac #1#2{{ #1\over #2 }}
\def \mm #1#2#3#4{\pmatrix #1 & #2 \cr #3 & #4 \endpmatrix}
\def\intl{\int\limits}

\def\({\left(}
\def\){\right)}

\def \C {{\Bbb C}}
\def \R {{\Bbb R}}
\def \Z {{\Bbb Z}}
\def \H {{\Cal H}}

\magnification=\magstep1

\NoBlackBoxes

\head
1.  Introduction and statement of results
\endhead

By a `converse theorem' we mean a uniqueness and existence statement
about a class of Dirichlet series.  A typical converse theorem
asserts that the only Dirichlet series with a given list of
properties are among those which have already been discovered.
The first converse theorem, proven by Hamburger \cite{H}
in 1922, states that the Riemann $\zeta$-function is
characterized by its functional equation.

Hecke showed that the $L$-functions associated with holomorphic
modular forms of even integral weight for the full modular group
satisfy certain functional equations, and conversely, the only Dirichlet
series satisfying these functional equations are $L$-functions
associated with modular forms.  This is the source of the term
`converse theorem.' 
We give some notation and then describe Hecke's result.

Throughout the paper we let
$$
F(s)=\sum_{n=1}^\infty\frac{a_n}{n^s}.
$$
We assume that $F(s)$ converges in some right half-plane and
continues to an entire function such that $\Gamma(s)F(s)$ is
entire and bounded in vertical strips.
This condition is denoted EBV.  We say that $F$ satisfies
a functional equation of degree~$2$, level~$N$, and
weight~$k$, if  
$$\eqalign{ 
\Phi(s) &= \left(\frac {\sqrt{N}}{2\pi} \right)^s
		\Gamma(s)F(s) \cr
&=\pm (-1)^{k/2}
{\Phi(k-s)}.
}$$

If $g$ is a function on the complex upper half-plane 
$\H = \{z\in \C \ :\ y>0\}$,  and 
$\gamma=\mm abcd \in GL_2(\R)^+ $, 
then we define the function
$g|_k\gamma$ by
$$
\left(g|_k\gamma\right) (z) = (\det\gamma)^{k/2} (cz+d)^{-k}
g\left({az+b\over cz+d}\right).
$$

The Hecke congruence group of level $N$ is defined by
$$
\Gamma_0(N) = \left\{ \mm abcd \in SL_2(\Z)\  :\  N|c\right\}.
$$
A function  $g:{\H}\to \C$ is called a cusp form of weight $k$
for $\Gamma_0(N)$ if 
$g|_k\gamma = g$ for all $\gamma\in \Gamma_0(N)$, and $g$ vanishes at
all cusps of $\Gamma_0(N)$.  
The space of cusp forms of weight $k$ is denoted $S_k(\Gamma_0(N))$.

Associated to the Dirichlet series $F(s)$ is a function on ${\H}$,
$$f(z)=\sum_{n=1}^\infty a_ne(nz).$$
Hecke's converse theorem relates properties of $F(s)$ to
properties of $f(z)$.

\proclaim{Hecke's converse theorem}  Suppose $N=1,2,3$ or $4$. 
If $F(s)$ is EBV and satisfies
a functional equation of degree~$2$, level~$N$, and weight~$k$, then
$f|_k\gamma = f$ for all $\gamma\in \Gamma_0(N)$.
\endproclaim

In fact it holds that $f\in S_k(\Gamma_0(N))$, but for now we are 
only concerned with the transformation properties of $f(z)$.

\demo {Proof}  By the Mellin inversion formula, if $y>0$,
$$\eqalign{ f(iy)&=\frac{1}{2\pi i} \intl_{(c)}
N^{-s/2}\Phi(s)y^{-s}~ds\cr
&=\frac{1}{2\pi i} \intl_{(c)}N^{(s-k)/2}\Phi(k-s)y^{s-k}~ds,}
$$

\noindent for $c>0$.  
The functional equation $\Phi(s)=\pm (-1)^{k/2}\Phi(k-s)$ gives
the transformation rule
$f(iy)=\pm N^{k/2}(iNy)^{-k}f(-1/iNy)$.  Since $f$ is holomorphic, this holds
for $y$ with positive real part.  In other words, $f|_kH_N = \pm f$,
where
$$
H_N = \mm {}{-1}N{} .
$$
Since $f(z) = f(z+1)$,  $f$ is invariant under
$$
\mm 11{}1  
\ \ \ \ \ \ \ \ \ \ \ \ \ \hbox{\rm and}
\ \ \ \ \ \ \ \ \ \ \ \ \
H_N \mm 1{-1}{}1 H_N^{-1} = \mm 1{}N1 .
$$
If $N=1,2,3$, or 4, those two matrices generate $\Gamma_0(N)$, 
proving the theorem.
\enddemo

If $N\ge 5$ then the above argument fails, and in fact the space
of functions satisfying the given conditions is infinite dimensional.
In order to get the desired conclusion that~$f$ is invariant under
$\Gamma_0(N)$, we must put further restrictions on~$F$.  Weil~\cite{W} 
conceived of the important idea of requiring that the twists
of~$F$ by Dirichlet characters satisfy an appropriate
functional equation.  
Later versions by Razar~\cite{Raz} and Li~\cite{Li}
reduced the number of twists to a finite number
depending on~$N$.  All subsequent converse theorems
for higher-rank groups \cite{JPS} \cite{P-S}
are built on the idea of requiring a functional 
equation for~$F$ and also for various twists of~$F$.

In this paper we have partial success at replacing the assumption
on twists of $F$ by the assumption of $F$ having an
Euler product of the appropriate form.  We say that $F$ has
an Euler product of degree~$2$, level~$N$, and weight~$k$, if
$$
F(s) = \prod_{p\ prime} F_p(s),
$$
where
$$ \eqalignno{
F_p(s) &= (1-a_p p^{-s} + p^{k-1-2s})^{-1}
& \hbox{if}\ \ p\nmid N\phantom{.}\ \ \ \ \ \cr
\cr
F_q(s) &= (1- q^{{k\over 2}-1-s})^{-1}
& \hbox{if}\ \ \ q\|N\phantom{.}\ \ \ \ \ \cr
\cr
F_q(s) &= 1
& \hbox{if}\ \ q^2|N.\ \ \ \ \ \cr
}$$

Our result is:

\proclaim{Theorem 1}  Let $5\le N\le 12$, or $14\le N\le 17$, or $N=23$, and
suppose $F(s)$ is EBV and has both a functional
equation and an Euler product of degree~$2$, level~$N$, and
weight~$k$. Then 
$f|_k\gamma=f$ for all $\gamma\in \Gamma_0(N)$.
\endproclaim

The proof only makes use of the Euler product
at a finite number of places, depending on $N$.

\proclaim{Corollary 1} Under the conditions of Theorem 1, 
$f\in S_k(\Gamma_0(N))$.
\endproclaim

The deduction that $f$ vanishes at the cusps of $\Gamma_0(N)$,
that is, $f\in S_k(\Gamma_0(N))$, is 
described in Section~5.

This work is motivated by the Selberg class of
Dirichlet series \cite{Se}.
One is led to believe that admission to this 
class is reserved for very special $L$-functions, and
probably all are associated with automorphic forms.
For instance,  Conrey and Ghosh~\cite{CG} show that the only elements
of the Selberg class of degree~$1$ are the Riemann
zeta-function and the Dirichlet
$L$-functions associated with primitive Dirichlet
characters.
One would like to show that the only elements of
the Selberg class of degree~$2$ are the $L$-functions
associated to cusp forms, both holomorphic
and non-holomorphic, 
which are eigenvalues of the Hecke operators $T_p$
for $p\nmid N$ and of the Atkin-Lehner operators
$U_q$ for $q\mid N$. 
This would require using the
Euler product condition of the Selberg axioms in place of
the twists required by Weil's theorem.
Our result makes use of the Euler product, but requires it to be
of a special form.

Our results are stated in the case where $F(s)$ `looks like'
the Dirichlet series associated to
a holomorphic cusp form.  All of the methods work equally
well in the case that $F(s)$ `looks like' the Dirichlet series associated to
a $GL_2$ Maass form.  The only extra step is verifying that
the conclusion of Lemma~5 holds when
$f(z)$ is an eigenfunction of the hyperbolic Laplacian.
A proof of this is given by B\"ockle~\cite{Bo}.

The paper is organized as follows.  In Section~2 we provide more
background information and then derive results based on the 
shape of the local factor of $F(s)$ at $p=2$.
In Section~3 we present some additional
 general methods.
In Section~4 we present ad-hoc methods for which we haven't
found an appropriate generalization.
In Section~5 we use the shape of the Euler product of $F(s)$
to show that $f(z)$ vanishes at the cusps of~$\Gamma_0(N)$.

\head
2. The local factor at $p=2$ 
\endhead

Recall that for $p$ prime the Hecke operator $T_p$ is defined
by
$$
T_p = \mm p{}{}1 + \sum_{a=0}^{p-1} \mm 1a{}p ,
$$
and for $q$ prime the Atkin-Lehner operator $U_q$ is
defined by
$$
U_q = \sum_{a=0}^{q-1} \mm qa{}q .
$$
We also put
$$
H_N = \mm {}{-1}N{} ,
\ \ \ \ \ \ \ 
P= \mm 11{}1  ,
\ \ \ \ \ \ \ 
W_N = \mm 1{}N1 .
$$

If a cusp form of weight $k$ on $\Gamma_0(N)$ is an eigenfunction of $H_N$, 
then its associated $L$-function will have a functional equation
of degree~$2$, level~$N$, and weight~$k$.  The reverse implication is 
also true, as we saw in the proof of Hecke's theorem.
If the cusp form is 
an eigenfunction of each $T_p$ for $p\nmid N$ and each $U_q$ for
$q|N$, then the $L$-function will have an Euler product of degree~$2$,
level~$N$, and weight~$k$.  It is easy to see that the reverse
holds also.  To summarize:

\proclaim {Lemma 1} 
If $F(s)$ has a functional equation of degree~$2$,
level~$N$, and weight~$k$, then
$f|_kP=f$ and $f|_kH_N=\pm f$, and so $f|_k W_N=f$.
If $F(s)$ has an Euler product of degree~$2$,
level~$N$, and weight~$k$, then
$$\eqalignno{
f|_k T_p &= a_p f & \hbox{\rm if}\ \ p\nmid N\phantom{.}\ \ \ \ \ \cr
\cr
f|_k U_q &= f|_k \mm q{}{}1 & \hbox{\rm if}\ \ \ q\|N\phantom{.}\ \ \ \ \ \cr
\cr
f|_k U_q &= 0  & \hbox{\rm if}\ \ q^2|N.\ \ \ \ \ \cr
}$$
\endproclaim

All of the information about the Dirichlet series $F(s)$ has been
translated to equivalent information about the function $f(z)$.
We will use this to deduce
that $f(z)$ is invariant under $\Gamma_0(N)$.

It is convenient to introduce
$$
\Omega_f =\{\omega \in \C[GL_2(\R)^+]\ :\ f|_k\omega = 0\} .
$$
Note that $\Omega_f$ is a right ideal in the group ring $\C[GL_2(\R)^+]$.
The goal of showing that $f(z)$ is invariant under $\Gamma_0(N)$
can be rewritten as showing $\gamma\equiv 1 \bmod \Omega_f$
for all $\gamma\in \Gamma_0(N)$, or equivalently, for a set of
$\gamma$ which generate $\Gamma_0(N)$.  
From now on, all congruences 
are assumed to be $\bmod \ \Omega_f$.

For each value of $N$ mentioned in Theorem 1, we will exhibit
a set of matrices which generate $\Gamma_0(N)$.  These were found by
using the generators of $\Gamma(1)$ and the coset representatives of
$\Gamma_0(N)$ in $\Gamma(1)$ to find a (large) generating set which
was then reduced down to a manageable size.  Generators for
$\Gamma_0(N)$ are also given by Chuman~\cite{Ch}; see also
the preprint by Ingle, Moore, and Wichert~\cite{IMW}. 

The next three lemmas
demonstrate how to use the shape of the Euler product at the 
prime $p=2$ to produce 
additional matrices for which $f(z)$ is invariant.  
The lemmas naturally correspond to the three cases $2\nmid N$,
$2\| N$, and $4|N$.

A useful calculation which will be used repeatedly is:
$$
H_N \mm ab{cN}d H_N \equiv \mm d{-c}{-bN}a .
$$
Let
$$
M_2 = \mm 21N{(N+1)/2} .
$$

\proclaim {Lemma 2}   If $P\equiv 1$, $H_N\equiv \pm 1$, and $T_2\equiv a_2$
for some $a_2\in \C$,  then $M_2\equiv 1$.
\endproclaim

\demo{Proof}  We are given
$$
\mm 2{}{}1 + \mm 1{}{}2 + \mm 11{}2 \equiv a_2 .
$$
Left multiplying and right multiplying by $H_N$ gives
$$
\mm 1{}{}2 + \mm 2{}{}1 + \mm 2{}N1 \equiv a_2 .
$$
This new congruence is valid because $H_N\equiv \pm 1$, and
$\Omega_f$ is a right ideal.  Subtract the two congruences to get
$$
\mm 2{}N1 
\equiv 
\mm 11{}2 . 
$$
Right multiply by $\mm 21{}1 $  and use $P\equiv 1$ to get $M_2\equiv 1$.
\enddemo

\demo{Proof of Theorem 1 for $N=5$, $7$, and $9$}
For those values of $N$, the group $\Gamma_0(N)$ has generators
$$
\Gamma_0(N)= \left\langle P,\ W_N,\ M_2 \right\rangle .
$$
By Lemmas 1 and 2, $f(z)$ is invariant under each of those matrices.
\enddemo

\proclaim{Lemma 3}  If $H_N\equiv \pm 1$ and $U_2\equiv \mm 2{}{}1$,
then
$$
\mm {-2}1N{-(N+2)/2} \equiv -1.
$$
\endproclaim

\demo{Proof}  We are given
$$
\mm 1{}{}1 + \mm 21{}2 \equiv \mm 2{}{}1 .
$$
Multiply on the right by $H_{2N}$ and use the relation
$$\eqalign{
H_{2N} &= H_N \mm 2{}{}1 \cr
&\equiv \pm \mm 2{}{}1 ,
}$$
to get
$$
\pm \mm 2{}{}1 + \mm 21{}2 H_{2N} \equiv \pm \mm 1{}{}1 .
$$
Combine the two congruences to get
$$
\mm 21{}2 H_{2N} \equiv \mp \mm 21{}2 .
$$
Right multiply by $\mm 21{}2 ^{-1}$ and left multiply by $H_N$ 
to get the stated relation. 
\enddemo

\demo{Proof of Theorem 1 for $N=6$ and $10$}
Let $A$ denote the matrix in Lemma~3, so $A\equiv -1 $. We have the
following lists of generators:
$$
\Gamma_0(6) = \left\langle P, \ W_6, \ A^{-1}W_6 A \right\rangle ,
$$
$$
\Gamma_0(10) = \left\langle P, \ W_{10}, \ (W_{10} A)^2 , \ 
H_{10} (W_{10}A)^2 H_{10}, \ A^{-1} W_{10}^{-1} A P^{-1} \right\rangle ,
$$
By Lemmas 1 and 3, $f(z)$ is invariant under those matrices.
\enddemo

\proclaim{Lemma 4} If $\,U_2\equiv 0$ then 
$$
\mm 21{}2 \equiv -1 .
$$
\endproclaim

\demo{Proof} Trivial.
\enddemo

\demo{Proof of Theorem 1 for $N=8$,  $12$ and $16$}
Let $B$ denote the matrix in Lemma 4,  so $B\equiv -1$.  We have 
the following lists of generators:
$$
\Gamma_0(8) = \left \langle  P,\ W_8, \ B^{-1}W_8 B  \right \rangle ,
$$
$$
\Gamma_0(12) = \left \langle P,\ W_{12}, \ B W_{12}^{-1} B , \    
H_{12} B^{-1} W_{12} B^{-1} H_{12}, \  B H_{12}B W_{12}^{-1} B H_{12} B
\right \rangle ,
$$
$$
\Gamma_0(16) = \left \langle  P, \ W_{16}, \ B W_{16}^{-1} B , \ 
(B H_{16})^4, \ (B^{-1} H_{16})^4  \right \rangle .
$$
By Lemmas 1 and 4, $f(z)$ is invariant under each of those matrices.
\enddemo

We have seen that for each $N$, the local factor of $F(s)$ at $p=2$
can be used to deduce invariance properties of $f(z)$ which are not
obtainable from Hecke's method.  For certain small $N$, this is
sufficient to deduce the invariance of $f(z)$ under all of
$\Gamma_0(N)$.  In the next two sections we make use of the local
factors at other primes to deduce further invariance properties of
$f(z)$.

\head 
3. General methods 
\endhead

We begin with a generalization of Lemma 2.

Let
$$
R_n={\sum_{1\le a \le n}}^{\!\!\!\prime}\ \mm na0n ,
$$
where ${\sum}^\prime$ means that the sum is over $(a,n)=1$.
Note that
$$T_p \mm p001 =R_p+\mm {p^2}{}{}1 + \mm 1{}{}1 .$$

\proclaim{Theorem 2} If $P\equiv 1$, and 
for each $p\mid n$
we have  $T_p\equiv \alpha $ for some $\alpha \in \C$,
 then
$H_NR_nH_{n^2N}\equiv R_n .$
\endproclaim

\demo{Proof} First consider the case where
$n=p^\lambda$. Using the fundamental identity 
$T_{p^\lambda}T_p\equiv T_{p^{\lambda+1}}+p^{k-1}T_{p^{\lambda-1}}$,
it is not difficult to prove by induction that
$$\eqalign{T_{p^\lambda}&=\sum_{j=0}^\lambda\sum_{b=0}^{p^j-1}
\mm {p^{\lambda-j}}b0{p^j}  \cr
&\equiv \alpha_\lambda ,
}$$
for some $\alpha_\lambda \in \C$. 
Multiplying  $T_{p^{\lambda-1}}\equiv \alpha_{\lambda-1}$ by
$\mm p{}{}1 $  gives
$$\sum_{j=0}^{\lambda-1}\sum_{b=0}^{p^j-1}
\mm {p^{\lambda-j}}b0{p^j}
\equiv \alpha_{\lambda-1} \mm p{}{}1 .$$
Subtracting this relation from
$T_{p^\lambda}\equiv \alpha_{\lambda}$, we obtain
$$T_{p^\lambda}-T_{p^{\lambda-1}}\mm p{}{}1
=\sum_{b=0}^{p^{\lambda-1}}\mm 1b{}{p^\lambda}
\equiv \alpha_\lambda -\alpha_{\lambda-1}\mm p{}{}1 . $$
Multiplying by $\mm p{}{}1 $
gives 
$$ \sum_{b=0}^{p^\lambda-1}\mm {p^\lambda}b{}{p^\lambda} 
\equiv \alpha_\lambda
\mm {p^\lambda}{}{}1 
-\alpha_{\lambda-1}\mm {p^{\lambda+1}}{}{}1 .$$
Now,
$$R_{p^\lambda}=\sum_{b=0}^{p^\lambda-1}
\mm {p^\lambda}b{}{p^\lambda} -
\sum_{b=0}^{p^{\lambda-1}-1}
\mm {p^{\lambda-1}}b{}{p^{\lambda-1}}  .$$
Therefore, by use of the previous relation twice,
$$R_{p^\lambda}\equiv 
\alpha_{\lambda}\mm {p^{\lambda}}{}{}1 -
\alpha_{\lambda-1}\mm {p^{\lambda+1}}{}{}1 -
\alpha_{\lambda-1}\mm {p^{\lambda-1}}{}{}1 +
\alpha_{\lambda-2}\mm {p^{\lambda}}{}{}1  .
$$ 
Now,
$$H_N \mm x{}{}y H_{Np^{2\lambda}} \equiv \mm {p^{2\lambda} y}{}{}x .$$ 
From this, it is easy to see that the right side of the
 above is invariant under left multiplication by
$H_N$ and right multiplication by $H_{Np^{2\lambda}}$.
This concludes the proof in the case that $n$ is a prime power.

To handle the general case, just note that if $(m_1,m_2)=1$, then
$$ 
R_{m_1m_2}\equiv R_{m_1}R_{m_2},
$$
and 
$$
\mm {m_1}{}{}1 R_{m_2} \equiv  R_{m_2} \mm {m_1}{}{}1 .
$$
\enddemo

Now we put Theorem~2 into a more usable form.
Let
$$\beta(x)=\mm 1x{}1 . $$

\proclaim{Corollary 2} 
If $P\equiv 1$, $H_N\equiv \pm 1$,  and 
for each $p\mid n$,
 $T_p\equiv \alpha $ for some $\alpha \in \C$,
then
$$
{\sum_{b(m)}}'\left(1-\mm m{-b}{-Nc}n 
\right)\beta\left(\frac bm
\right)\equiv 0 ,
$$
where the sum is over a set of reduced residues $b$ modulo $m$, 
and where $c$ and $n$ are integers depending on $m$, $b$ and $N$ such that
$mn-bcN=1$.
\endproclaim

\demo{Proof} Let
$$\gamma(b,c)=\mm m{-b}{-Nc}n \in\Gamma_0(N).$$
It is an easy calculation to check that
$$\eqalign{ \beta\left(\frac cm \right)H_{Nm^2}&=
H_N\gamma(b,c)\beta\left(\frac bm
\right)\mm N{}{}N \cr
&\equiv \pm \gamma(b,c)\beta\left(\frac bm \right)  .
}$$
Now,
$$
R_m\equiv{\sum_{c(m)}}'\beta\left(\frac cm \right).
$$
Therefore, by Theorem 2,
$$\eqalign{ {\sum_{c(m)}}'\beta\left(\frac cm \right)
&\equiv \pm  {\sum_{c(m)}}'\beta\left(\frac cm \right)H_{Nm^2}\cr 
&\equiv {\sum_{c(m)}}'\gamma(c,b)\beta \left(\frac cm \right) .
} $$
This relation implies Corollary 2.
\enddemo

If $m=2$ then Corollary 2 is equivalent to Lemma~2.

The reduction from Corollary 2 to invariance properties
of $f(z)$ uses ideas from the proof of Weil's
converse theorem as described in Ogg's book \cite{Ogg}. We
quote Proposition 3 from that book for convenience:

\proclaim{Lemma 5} Suppose $f$ is holomorphic in $\H$ and
$\varepsilon\in GL_2(\R)^+$ is elliptic.  If \hbox{$f|_k\varepsilon=f$},
then either $\varepsilon$ has finite order, or $f$ is constant.
\endproclaim

The typical way we apply Lemma 5 is to use Corollary 2 to
first prove 
$$
1-\gamma\equiv (1-\gamma)\varepsilon
$$
for some $\gamma\in SL_2(\Z)$ and some elliptic $\varepsilon \in
GL_2(\R)^+$ which is not of finite order. We then conclude 
$\gamma\equiv 1$.

We illustrate the method in the case $N=11$. By Corollary~2,
we have
$$\left(1-\mm 3{-1}{-11}4 \right)\beta(1/3)
+\left(1- \mm 31{11}4 \right)\beta(-1/3)
\equiv 0 $$
and
$$\left(1- \mm 4{-1}{-11}3 \right)\beta(1/4)
+\left(1- \mm 41{11}3 \right)\beta(-1/4)
\equiv 0. $$
Therefore,
$$\eqalign{
1- \mm 3{-1}{-11}4 
&\equiv-\left(1- \mm 31{11}4 \right)\beta\left(-\frac 23 \right)\cr 
&=\left(1- \mm 4{-1}{-11}3 \right) \mm 31{11}4 \beta\left(-\frac 23 \right)\cr 
&\equiv -\left(1- \mm 41{11}3 \right)\beta\left(-\frac{2}{4} \right)
\mm 31{11}4 \beta\left(-\frac 23 \right)\cr 
&=\left(1-\mm 3{-1}{-11}4 \right)\mm 41{11}3 \beta\left(
-\frac 24 \right)\mm  31{11}4 \beta\left(-\frac 23 \right) .
}$$
However,
$$
\mm 41{11}3
\beta\left(-\frac{2}{4} \right)
\mm 31{11}4 \beta\left(-\frac 23 \right)
= \mm 1{-2/3}{11/2}{-8/3} 
$$
is elliptic but not of finite order.
So,
$$ 
\mm 3{-1}{-11}4 \equiv 1 .
$$

We will use the above calculation to prove Theorem 1 for $N=11$, and
then we will describe a generalization of the method.

Let
$$
M_{m,b}=M_{m,b}(N) = \mm mb{cN}d ,
$$
where $0<2|c|<|m|$.  Also put $M_m=M_{m,1}$.  
This is ambiguous only if $m=\pm2$, in which case we use our
previous definition of $M_2$, and put $M_{-2}= W_N^{-1}M_2 P^{-1}$.

\demo{Proof of Theorem 1 for $N=11$} We have generators
$$
\Gamma_0(11) = \left\langle P,\ M_2,\ M_3 \right\rangle .
$$
By Lemmas 1 and 2, and Corollary 2 and the above calculation, 
$f(z)$ is invariant under
the given matrices.
\enddemo

We summarize the above argument as follows. First of all, for 
$$\gamma=\mm abcd $$
let
$$\gamma'=\mm a{-b}{-c}d .$$
Clearly $\gamma$ and $\gamma'$ have the same determinant and
belong to the same group $\Gamma_0(N)$.

Now let $m$ and $n$ be two of the integers 3, 4, and 6. (These are
the integers  with $\varphi(m)=\varphi(n)=2$.)
Let
$$
\gamma=\mm m1Nn \in \Gamma_0(N).
$$
Thus, $N$ will be one of 8, 11, 15, 17, 23, or 35. 
Corollary~2 implies 
$$
(1-\gamma)\beta(-1/m)+(1-\gamma')\beta(1/m)\equiv 0,
$$
and also 
$$(1-\gamma'{}^{-1})\beta(-1/n)+(1-\gamma^{-1})\beta(1/n)\equiv 0.$$
Therefore,
$$\eqalign{
1-\gamma &\equiv -(1-\gamma')\beta(2/m)\cr 
&=(1-\gamma'{}^{-1})\gamma'\beta(2/m)\cr 
&\equiv -(1-\gamma^{-1})\beta(2/n)\gamma'\beta(2/m)\cr 
&=(1-\gamma)\gamma^{-1}\beta(2/n)\gamma'\beta(2/m).
}$$
But
$$
\gamma^{-1}\beta(2/n)\gamma'\beta(2/m)=
		\mm 1{\frac2m}{-\frac{2N}n}{\frac4{mn}-3} 
$$
is elliptic of infinite order.
Hence,
$$\gamma \equiv 1, $$
and similarly for $\gamma'$.

What we have shown is:

\proclaim{Corollary 3} Under the conditions of Corollary $2$, 
if 
$\,\varphi(m)=\varphi((N+1)/m)=2$, then
$M_m\equiv 1$ and $\,M_{-m}\equiv 1$.
\endproclaim

\demo{Proof of Theorem 1 for $N=17$}  We have generators, 
$$
\Gamma_0(17) = \left\langle P,\ W_{17},\ M_2,\ M_3,\ M_6 \right\rangle .
$$
By Lemmas 1 and 2, and Corollary 3, $f(z)$ is invariant under
the above matrices.
\enddemo

\head
4. Ad-hoc methods
\endhead

The previous two sections gave general methods for finding
matrices $\gamma\in \Gamma_0(N)$ such that $f|_k\gamma=f$. This
was sufficient to generate $\Gamma_0(N)$ for a few $N$.
In this section we use ad-hoc methods to deduce the invariance
of $f(z)$ for various other matrices.  It may be that discovering a 
general scheme behind these seemingly ad-hoc methods 
could lead to a proof of the general case.

The main weakness in our method of using Corollary~2 is its intractability
when $\varphi(m)>2$.  The next two proofs exhibit a `bootstrap'
feature, where expressions with $\varphi(m)>2$ are first reduced
down to simpler expressions, and then these simpler expressions
are used in a manner similar to the proof of Corollary~3.

\demo{Proof of Theorem 1 for $N=14$}  We have generators,
$$
\Gamma_0(14)=\left\langle P,\ W_{14},\ M_{3},\ M_{-3},\ M_{13,6} 
\right\rangle  .
$$
By Lemma 1, $f$ is invariant under $P$ and $W_{14}$.  Writing
the conclusion of Lemma~3 as $A\equiv -1$, we have that $f$ is
invariant under
$$
A^{-1}W_{14}^{-1} A  P = M_{13,6} .
$$
It remains to show invariance under the two other generators.

Corollary~2, with $m=3$, gives
$$\eqalign{
\(1-\mm 3{-1}{-14}5 \)\beta\({1\over 3}\) +
	\(1-\mm 31{14}5 \) \beta\(-{1\over 3}\) &\equiv 0, \cr
}$$
which we write as
$$\eqalign{
	\(1-\gamma_1\)\beta\({1\over 3}\)
			+ \(1-\gamma_2\)\beta\(-{1\over 3}\) &\equiv 0,
}$$
temporarily putting $\gamma_1=M_{-3}$ and $\gamma_2=M_3$.
Corollary~2, with $m=5$, gives
$$\eqalign{
	&\(1-\mm 5{-2}{28}{-11} \)\beta\({2\over 5}\) +
		\(1-\mm 5{-1}{-14}3 \)\beta\({1\over 5}\) \cr
\cr
&\phantom{XXXXXX}+
	\(1-\mm 51{14}3 \)  \beta\(-{1\over 5}\)
	+ \(1-\mm 52{-28}{-11} \) \beta\(-{2\over 5}\) \equiv 0.
}$$
We can reduce the $m=5$ expression by noting the following:
$$\leqalignno{
(AP)^2 &= \mm 52{-28}{-11} \cr
	&\equiv 1,\cr
&&\hbox{and}\cr
(W_{14}A)^2 &= \mm 5{-2}{28}{-11} \cr
	&\equiv 1.
}$$
So the $m=5$ expression implies
$$
\(1-\gamma_2^{-1}\) \beta\({1\over 5}\) +
\(1-\gamma_1^{-1}\) \beta\(-{1\over 5}\)
\equiv 0.
$$
Now proceed exactly as in  the proof of Corollary~3:
$$\eqalign{
1-\gamma_1&\equiv - \(1-\gamma_2\)\,\beta\(-{2\over 3}\)\cr
	&= \(1-\gamma_2^{-1}\)\,\gamma_2\,\beta\(-{2\over 3}\)\cr
	&\equiv - \(1-\gamma_1^{-1}\)\,\beta\(-{2\over 5}\)\,
	\gamma_2\,\beta\(-{2\over 3}\)\cr
	&=\(1-\gamma_1\)\,\gamma_1^{-1}\,\beta\(-{2\over 5}\)\,
		\gamma_2\,\beta\(-{2\over 3}\)\cr
	&=\(1-\gamma_1\) \mm 1{-2/3}{28/5}{-41/15}  .
}$$
That last matrix is elliptic of infinite order, so $\gamma_1\equiv 1$.
Therefore, $\gamma_2\equiv 1$.
This finishes the proof of Theorem~1 for $N=14$.  
\enddemo

Note that for the above proof we used the local factors of $F(s)$ at
$p=2$, 3, and~5, while in all other cases we only used the factors
of $F(s)$ at $p=2$ and~3.

\demo{Proof of Theorem 1 for $N=15$} We have generators,
$$
\Gamma_0(15)=
\left\langle
P,\ W_{15},\ M_2,\ M_4,\ M_{11,4}
\right\rangle .
$$
From Lemmas 1 and 2, and Corollary 3, we have that $f(z)$ is
invariant under the first four generators.  It remains
to prove invariance under $M_{11,4}$.

Corollary 2, with $m=8$, gives
$$\eqalign{
&\(1 - \mm 8{-3}{-45}{17} \)  \beta\({3\over 8}\) +
\(1 - \mm  8{-1}{-15}2 \) \beta\({1\over 8}\) \cr
&\phantom{XXXXXXXX}+
\(1 - \mm 81{15}2  \)  \beta\(-{1\over 8}\) +
\(1 - \mm 83{45}{17} \)  \beta\(- {3\over 8}\) \equiv 0.
}$$
We can reduce this expression by noting that
$$\leqalignno{
M_2^{-1}&=\mm 8{-1}{-15}2 \cr
	&\equiv 1, \cr
&&\hbox{and}\cr
PM_2^{-1}W_{15} &= \mm 81{15}2  \cr
	&\equiv 1.
}$$
And we also have
$$\leqalignno{
P M_2^{-1} \mm 83{45}{17} &= M_{11,4}\cr
&&\hbox{and}\cr
M_2^{-1} W_{15} \mm 8{-3}{-45}{17} &=  M_{11,4}^{-1}  .
}$$
So the $m=8$ relation gives
$$\leqalignno{
1-M_{11,4} &\equiv - (1-M_{11,4}^{-1}) \beta\(\frac 34\) \cr
&= (1-M_{11,4}) M_{11,4}^{-1}\, \beta\(\frac 34\) \cr
&= (1-M_{11,4}) \mm {11}{17/4}{-30}{-23/2}  .
}$$
The last matrix is elliptic of infinite order, so $M_{11,4}\equiv 1$.
This finishes the proof of Theorem 1 for  $N=15$.
\enddemo

\demo{Proof of Theorem 1 for $N=23$} We have generators,
$$
\Gamma_0(23) =
\left\langle
P,\ W_{23},\ M_2,\ M_3,\ M_4,\ M_6
\right\rangle .
$$
By Lemmas 1 and 2, and Corollary 3, we obtain the invariance of $f(z)$
under each of the above matrices except $M_3$.

By Corollary 2, with $m=3$, we get
$$
\(1-\mm  31{23}8 \)\beta \(\frac {-1}3 \) + 
\(1-\mm 3{-1}{23}{-8} \) \beta \(\frac 13 \) \equiv 0 ,
$$
which we can rewrite as
$$
(1-M_3)\beta \(\frac {-1}3 \) + 
(1-M_{-3}) \beta \(\frac 13 \) \equiv 0.
$$
Let
$$
\gamma = \mm {10}3{23}7 .
$$
One can check that 
$W_{23}M_{-2}\gamma=M_3$ and $M_{-6}M_6M_{-3}\gamma=I$.
In particular, $\gamma\equiv M_3$ and $M_{-3}\gamma\equiv 1$. Consequently, 
$$
1-M_3\equiv -(1-M_{-3})\gamma .
$$
Combining the two expressions gives
$$ 
1-M_3\equiv (1-M_3)\beta(-2/3)\gamma .
$$
But
$$\beta(-2/3)\gamma=\mm {-16/3}{-5/3}{23}7 $$
is elliptic of infinite order. So, $M_3\equiv 1$ as desired. 
This finishes the proof of Theorem~1 for $N=23$.
\enddemo

This finishes the proof of Theorem~1.  It remains to show that
$f(z)$ is actually a cusp form.

\head
5. Vanishing at the cusps
\endhead

We have assumed that $F(s)$ converges in some right half-plane.
It follows from this, and the invariance of $f(z)$ under $\Gamma_0(N)$,
that $f(z)$ is holomorphic at the cusps of $\Gamma_0(N)$.
If one assumes further that $F(s)$ converges for $\sigma>k-\delta$, for
some $\delta>0$, then it follows that
$f(z)$ vanishes at the cusps of $\Gamma_0(N)$.
See Proposition~1 and the Lemma on page V-14 of Ogg's book~\cite{Ogg}.
Note that the analog of the Ramanujan conjecture gives bounds
on the $a_n$ which are much stronger than is required for this method.
We will show that no such assumption is needed.  The vanishing of
$f(z)$ at the cusps of $\Gamma_0(N)$ follows from the previously
assumed functional equation and Euler product of $F(s)$.

The following, in combination with Lemma~1, gives Corollary~1
as a consequence of Theorem~1.

\proclaim{Theorem 3} Let $N=2^kN'$, with $0\le k\le3$ and $N'$ odd
and squarefree, and put $\displaystyle{f(z)=\sum_{n=1}^\infty a_ne(nz)}$.  
If $f|_k\gamma = f$ for all 
$\gamma\in\Gamma_0(N)$, and $f|_kH_N=\pm f$, and
$$\eqalignno{
f|_k U_q &= f|_k \mm q{}{}1 & \hbox{\rm if}\ \ \ q\|N\phantom{.}\ \ \ \ \ \cr
\cr
f|_k U_q &= 0  & \hbox{\rm if}\ \ q^2|N,\ \ \ \ \ \cr
}$$
then $f(z)$ vanishes at the cusps of $\Gamma_0(N)$.
\endproclaim

The values of $N$ in Theorem~3 are exactly those for which 
$\{1/r\ :\ r|N\}$ is a set of cusps for $\Gamma_0(N)$.  Note 
that in this case the cusp $1/r$ has width $N/r$.

The following Lemma follows from the Chinese remainder theorem.
It is implicit in the classification of the cusps of $\Gamma_0(N)$.

\proclaim{Lemma 6} If $r|N$ and $(ar,N/r)=1$, then there exists
$\gamma\in\Gamma_0(N)$ such that
$$
\gamma \mm 1{}{ar}1 = \mm 1b{-r}d 
$$
for some $b,d\in \Z$.
\endproclaim

\demo{Proof of Theorem 3}  Suppose $r|N$.  We show that $f(z)$
vanishes at the cusp $1/r$.

Put $N=rQ$.  By taking linear combinations of $U_q$ for $q|Q$, we can
obtain, for some $\beta_q\in \Z$,
$$
{\sum_{1\le a\le Q}}^{\!\!\!\prime} \mm 1{a/Q}{}1 \equiv \sum_{q|Q} \beta_q \mm q{}{}1 .
$$
Left multiply and right multiply by $H_N$ to get
$$
{\sum_{1\le a\le Q}}^{\!\!\!\prime} \mm 1{}{ar}1 \equiv \sum_{q|Q} \beta_q \mm 1{}{}q .
$$
So by Lemma~6 we have
$$
{\sum_{1\le a\le Q}}^{\!\!\!\prime} \mm 1{b_a}{-r}{d_a} 
	\equiv \sum_{q|Q}  \beta_q\mm 1{}{}q ,
$$
which can be rewritten as
$$
\mm 1{}{-r}1 {\sum_{1\le a\le Q}}^{\!\!\!\prime} \mm 1{b_a}{}1 
\equiv \sum_{q|Q}  \beta_q\mm 1{}{}q .
$$
By definition, that congruence says,
$$\eqalign{
\(f\Big|_k \mm 1{}{-r}1 {\sum_{1\le a\le Q}}^{\!\!\!\prime} 
			\mm 1{b_a}{}1 \) (z)
&= \(f \Big|_k \sum_{q|Q}  \beta_q\mm 1{}{}q \) (z) \cr
&= \sum_{n=1}^\infty c_n e(nz/Q) ,
}$$
for some $c_n\in \C$.
Now, if $\alpha_0$ is the constant term in the expansion
of $f|_k \mm 1{}{-r}1 $, then the constant term on the left side of
the above expression is $\varphi(Q)\alpha_0$.  The expression on the right
side above has no constant term, so $\alpha_0=0$.  In other words, $f(z)$
vanishes at the cusp $1/r$.  This proves Theorem~3.

\enddemo

The only values of $N$ in Theorem~1 which are not covered by
Theorem~3 are $N=9$ and~$16$.  Using the equality
$$\eqalign{
H_{M^2} \mm 1{}{M}1 &= \mm {-M}{-1}{M^2}{} \cr
&\equiv \mm 1{}{-M}1 \mm 1{1/M}{}1 ,
}$$
it is possible to modify the above proof to give the 
conclusion of Theorem~3  for $N=9$ and $16$ 
in the particular case $f|_k H_N =f$, that is, when $F(s)$ has
sign $+1$ in its functional equation.  If $F(s)$ has sign $-1$ 
in its functional equation, then the above method fails and we
must make an assumption on the growth of the coefficients of
$F(s)$ to obtain Corollary~1.  This extra assumption cannot
be eliminated, as demonstrated by the following example. 
The Dirichlet series $L(s,\chi_3)L(s-1,\chi_3)$, where $\chi_3$ is 
the Dirichlet character mod~3, is entire and has a functional
equation and Euler product of degree~2, weight~2, and level~9.
Its functional equation has sign~$-1$, and it is associated
to an Eisenstein series of weight~2 on~$\Gamma_0(9)$.  Similar
examples exist for the other values of $N$ not covered by
Theorem~3.

\Refs

\item{[Bo]}  G.~B\"ockle, {\it An elementary proof of the analog
of Weil's converse theorem for Maass wave forms}, preprint.

\item{[CG]} J.B. Conrey and A. Ghosh, {\it On the Selberg class
of Dirichlet series: small degrees}, Duke Math.~J. {\bf 72}, (1993),
673-693.

\item{[Ch]} Y. Chuman,  {\it Generators and relations of
$\Gamma_0(N)$}, J. Math. Kyoto Univ.
{\bf 13}
(1973).

\item{[H]} H.~Hamburger, {\it \"Uber die Riemannsche Funktionalgleichung
der $\zeta$-Funktion}, Math. Zeit {\bf 13} (1922).

\item{[IMW]} Ingle, Moore, and Wichert, {\it Explicit generating sets for
$\Gamma_0(N)$, $N\le 100$}, preprint.

\item{[JPS]}
H. Jacquet, I.~Piatetski-Shapiro, and J.  Shalika,
{\it Automorphic forms on $GL(3)$ I,  II},
Ann.~Math. {\bf 109} (1979).

\item{[Li]} W.~Li, {\it On a theorem of Hecke-Weil-Jacquet-Langlands},
Proc. Symposium in Analytic Number Theory, Durham, England, 1979.

\item{[Ogg]} 
A. Ogg,
{\it  Modular Forms and Dirichlet Series},
W.A.~Benjamin, New York,  1969.

\item{[P-S]}  I. Piatetski-Shapiro,
{\it The converse theorem for $GL(n)$},
Festschrift in Honor of I.~Piatetski-Shapiro
vol. II, 1990, 185-196.

\item{[Raz]} 
M. Razar,
{\it Modular forms for $\Gamma_0(N)$ and Dirichlet series},
Trans. Amer. Math. Soc 
{\bf 231}
(1977), 
489-495.

\item{[Se]} A.  Selberg, {\it  Old and new results and conjectures about
a class of Dirichlet series}, Proceedings of the Amalfi Conference
on Analytic Number Theory 1989,
Univ. Salerno, Salerno (1992), 367-385.

\item{[W]} A. Weil,
{\it Uber die Bestimmung Dirichletscher Reihen durch
Funktionalgleichungen},
Math.~Ann. {\bf 168} (1967).

\endRefs

\enddocument